# On Hardy type spaces in strictly pseudoconvex domains and the density, in these spaces, of certain classes of singular functions

by

Kyranna Kioulafa

## Abstract

In this paper we prove generic results concerning Hardy spaces in one or several complex variables. More precisely, we show that the generic function in certain Hardy type spaces is totally unbounded and hence non-extentable, despite the fact that these functions have non-tangential limits at the boundary of the domain. We also consider local Hardy spaces and show that generically these functions do not belong – not even locally – to Hardy spaces of higher order. We work first in the case of the unit ball of $\mathbb{C}^n$ where the calculations are easier and the results are somehow better, and then we extend them to the case of strictly pseudoconvex domains.

## 1. Introduction and preliminaries

Hardy type spaces for domains in $\mathbb{C}^n$ and the boundary behaviour of the functions belonging to these spaces have been extensively studied by various authors. See for example [5, 7, 9, 11, 12, 13]. More relevant to the results of this paper are those in [3], where analogous theorems were proved with Bergman type spaces in place of Hardy spaces. We also refer to [3] and the references given there for the techniques used in this present paper.

Let $\mathbb{B} = \{z \in \mathbb{C}^n : |z| < 1\}$ be the unit ball in $\mathbb{C}^n$. Let us recall that the Hardy space $H^p(\mathbb{B})$, $1 \leq p < \infty$, is defined to be the set of holomorphic functions $f : \mathbb{B} \to \mathbb{C}$ such that

$$\|f\|_p = \sup_{r<1} \left( \int_{\zeta \in \partial \mathbb{B}} |f(r\zeta)|^p d\sigma(\zeta) \right)^{1/p} < +\infty,$$

where $d\sigma$ is the Euclidean surface area measure on the sphere $\partial \mathbb{B}$, and that with the norm $\|\cdot\|_p$, $H^p(\mathbb{B})$ becomes a Banach space. We also recall that if a sequence $f_m \in H^p(\mathbb{B})$ converges to $f$, in the above norm, then $f_m$ converges to $f$ also uniformly on compact subsets of $\mathbb{B}$. Indeed this follows from the inequality

$$\sup_{z \in K} |f(z)| \leq C(p, K) \|f\|_p,$$

with $K$ being a compact subset of $\mathbb{B}$ and $C(p, K)$ is a constant depending on $p$ and $K$. (See [9, Theorem 7.2.5].) Also $H^\infty(\mathbb{B})$ is the Banach space of bounded holomorphic functions $f : \mathbb{B} \to \mathbb{C}$, with the norm $\|f\|_\infty = \sup_{z \in \mathbb{B}} |f(z)|$.

For each $q > 1$, we also consider the space $\bigcap_{1 \leq p < q} H^p(\mathbb{B})$, which becomes a complete metric space with the metric

$$d(f,g) = \sum_{j=1}^{\infty} \frac{1}{2^j} \frac{\|f-g\|_{p_j}}{1+\|f-g\|_{p_j}},$$

where $1 < p_1 < p_2 < \cdots < p_j < \cdots < q$ and $p_j \to q$ ($j \to \infty$). Let us point out that the topology induced by this metric in the space $\bigcap_{1 \leq p < q} H^p(\mathbb{B})$ is independant of the choice of the sequence $p_j$. As a matter of fact, a sequence $f_k$ converges to $f$, in $\bigcap_{1 \leq p < q} H^p(\mathbb{B})$, if and only if $\|f_k - f\|_p \to 0$, for every $p < q$.

In the first part of the paper we will show that generically the functions in the space $\bigcap_{1 \leq p < q} H^p(\mathbb{B})$ are totally unbounded (Theorem 2.1). A function $f : \mathbb{B} \to \mathbb{C}$ is called totally unbounded if for every $\zeta \in \partial \mathbb{B}$ and every $\varepsilon > 0$, the restriction $f|_{B(\zeta,\varepsilon) \cap \mathbb{B}}$ of $f$ to the set $B(\zeta,\varepsilon) \cap \mathbb{B} = \{z \in \mathbb{B} : |z - \zeta| < \varepsilon\}$ is unbounded, i.e.,

$$\sup_{z \in B(\zeta,\varepsilon) \cap \mathbb{B}} |f(z)| = \infty.$$

Following a suggestion of Nestoridis, we also consider local Hardy spaces $H^p(\mathbb{B}, G)$, for open subsets $G$ of the sphere $\partial \mathbb{B}$ (the precise definition is given in section 2) as another way of measuring how singular a holomorphic function is, near a boundary point. One extreme of this concept is the property of being totally unbounded. (A similar definition – although slightly different – is given in [7].) In this context we show that generically the functions in the space $\bigcap_{1 \leq p < q} H^p(\mathbb{B})$ do not belong to any local $H^q$ – space (Theorem 2.4).

In sections 4 and 5, we will extend these results from the ball to the case of strictly pseudoconvex domains. In this more general case we have to modify the definition of local Hardy spaces which we give in the case of the ball, and consider the space $H^p(\Omega, U)$, where $U$ is an open subset of $\mathbb{C}^n$ so that $U \cap (\partial \Omega) \neq \varnothing$. (See section 3.)

## 2. The case of the unit ball of $\mathbb{C}^n$

In this section we will first prove the following theorem.

**Theorem 2.1.** *Let $q \in \mathbb{R} \cup \{\infty\}$, $q > 1$. Then the set of the functions in the space $\bigcap_{1 \leq p < q} H^p(\mathbb{B})$ which are totally unbounded in $\mathbb{B}$ is dense and $\mathcal{G}_\delta$ in this space.*

The proof of this theorem will be based on the following lemmas.

**Lemma 2.2.** *For each point $\zeta \in \mathbb{S} = \partial \mathbb{B}$, we consider the functions*

$$f_\zeta(z) = \frac{1}{1 - \langle z, \zeta \rangle} = \frac{1}{1 - \sum_{j=1}^n \overline{\zeta}_j z_j}, \quad h_\zeta(z) = \log f_\zeta(z) \text{ and } \varphi_{q,\zeta}(z) = \exp\left[\frac{n}{q} h_\zeta(z)\right] \text{ defined for } z \in \mathbb{B}.$$

*Then (i) $f_\zeta \in \bigcap_{1 \leq p < n} H^p(\mathbb{B})$ and $f_\zeta \notin H^n(\mathbb{B})$,*

*(ii) $h_\zeta \in \bigcap_{1 \leq p < \infty} H^p(\mathbb{B})$ and $h_\zeta \notin H^\infty(\mathbb{B})$,*

*(iii) $\varphi_{q,\zeta} \in \bigcap_{1 \leq p < q} H^p(\mathbb{B})$ and $\varphi_{q,\zeta} \notin H^q(\mathbb{B})$ for $1 < q < \infty$.*

**Proof.** By [9, Proposition 1.4.10], if $p < n$, the integral

$$\int_{w \in \mathbb{S}} \frac{d\sigma(w)}{|1 - \langle z, w \rangle|^p},$$

as a function of $z$, remains bounded for $z \in \mathbb{B}$, and, therefore since $r\zeta \in \mathbb{B}$ for $r < 1$,

$$\sup_{0 < r < 1} \int_{z \in \mathbb{S}} |f_\zeta(rz)|^p d\sigma(z) = \sup_{0 < r < 1} \int_{z \in \mathbb{S}} \frac{d\sigma(z)}{|1 - \langle rz, \zeta \rangle|^p} = \sup_{0 < r < 1} \int_{z \in \mathbb{S}} \frac{d\sigma(z)}{|1 - \langle r\zeta, z \rangle|^p} < \infty.$$

Thus $f_\zeta \in \bigcap_{1 \leq p < n} H^p(\mathbb{B})$. Next we show that

$$\sup_{0 < r < 1} \int_{z \in \mathbb{S}} |f_\zeta(rz)|^n d\sigma(z) = \sup_{0 < r < 1} \int_{z \in \mathbb{S}} \frac{d\sigma(z)}{|1 - \langle rz, \zeta \rangle|^n} = \infty.$$

Indeed, by [9, Proposition 1.4.10], the integral

$$\int_{z \in \mathbb{S}} \frac{d\sigma(z)}{|1 - \langle w, z \rangle|^n} \text{ behaves as } \log \frac{1}{1 - |w|^2} \text{ for } w \in \mathbb{B},$$

and therefore

$$\sup_{0 < r < 1} \int_{z \in \mathbb{S}} \frac{d\sigma(z)}{|1 - \langle rz, \zeta \rangle|^n} = \sup_{0 < r < 1} \int_{z \in \mathbb{S}} \frac{d\sigma(z)}{|1 - \langle r\zeta, z \rangle|^n} = \sup_{0 < r < 1} \log \frac{1}{1 - r^2} = \infty.$$

This proves (i). Next, observing that $\text{Re}(1 - \langle z, \zeta \rangle) > 0$, for $z \in \mathbb{B}$, we see that $\text{Re} f_\zeta(z) > 0$ and therefore we may define $h_\zeta(z) = \log f_\zeta(z)$ using the principal branch of the logarithm with $-\pi < \arg \leq \pi$. Then $|\text{Im}[\log f_\zeta(z)]| < \pi/2$, i.e., $h_\zeta(z) = \log |f_\zeta(z)| + i\theta(z)$ with $|\theta(z)| < \pi/2$. It follows that if the point $rz \in \mathbb{B}$ and is sufficiently close to $\zeta$,

$$|h_\zeta(rz)|^p = \left| \log \frac{1}{|1 - \langle rz, \zeta \rangle|} + i\theta(rz) \right|^p = \left[ \left( \log \frac{1}{|1 - \langle rz, \zeta \rangle|} \right)^2 + \theta^2(rz) \right]^{p/2} \preceq (k!)^{p/k} \frac{1}{|1 - \langle rz, \zeta \rangle|^{p/k}},$$

where we used the inequality $(\log x)^p \leq (k!)^{p/k} x^{p/k}$ which holds for $x > 1$, $p \geq 1$ and $k \in \mathbb{N}$. (We also used the fact that, since $|1 - \langle rz, \zeta \rangle| > 0$ for $rz \in \overline{\mathbb{B}}$ away from the point $\zeta$, the quantity $|\log|1 - \langle rz, \zeta \rangle||$ is bounded.) Fixing a $p < \infty$ and choosing $k > p/n$, we see (using also (i)) that $h_\zeta \in H^p(\mathbb{B})$ whence we obtain $h_\zeta \in \bigcap_{1 \leq p < \infty} H^p(\mathbb{B})$. Since obviously $\lim_{z \in \mathbb{B}, z \to \zeta} h_\zeta(z) = \infty$, (ii) follows.

Finally observing that $|\varphi_{q,\zeta}| = |f_\zeta|^{n/q}$, we easily obtain (iii).

The following lemma is proved in [10].

**Lemma 2.3.** *Let $\mathcal{V}$ be a topological vector space over $\mathbb{C}$, $X$ a non-empty set, and let $\mathbb{C}^X$ denote the vector space of all complex-valued functions on $X$. Suppose $T : \mathcal{V} \to \mathbb{C}^X$ is a sublinear operator with the property that, for every $x \in X$, the functional $T_x : \mathcal{V} \to \mathbb{C}$, defined by $T_x(f) := T(f)(x)$, for $f \in \mathcal{V}$, is continuous. Let $\mathcal{E} = \{f \in \mathcal{V} : T(f) \text{ is unbounded on } X\}$. Then either $\mathcal{E} = \varnothing$ or $\mathcal{E}$ is dense and $\mathcal{G}_\delta$ set in the space $\mathcal{V}$.*

***Proof of Theorem 2.1.*** Let us consider a 'small' ball $B$ whose center lies on $\partial \mathbb{B}$, and let us set $X = B \cap \mathbb{B}$ and $\mathcal{V} = \bigcap_{p < q} H^p(\mathbb{B})$. We define the linear operator

$$T : \mathcal{V} \to \mathbb{C}^X \text{ with } T(f)(z) = f(z) \text{ for } z \in X \text{ and } f \in \mathcal{V}.$$

For each fixed $z \in X$, the functional $T_z : \mathcal{V} \to \mathbb{C}$ defined by $T_z(f) = f(z)$, $f \in \mathcal{V}$, is continuous. It is easy to see that the set $\mathcal{E} = \{f \in \mathcal{V} : T(f) \text{ is unbounded on } X\}$ in this case is equal to

$$\mathcal{E}(B) = \left\{ f \in \bigcap_{1 \leq p < q} H^p(\mathbb{B}) : \sup_{z \in B \cap \mathbb{B}} |f(z)| = \infty \right\}.$$

Also, by Lemma 2.2(ii), $\mathcal{E}(B) \neq \emptyset$, since $h_\zeta \in \mathcal{E}(B)$ for $\zeta \in B \cap \partial \mathbb{B}$. Therefore, by Lemma 2.3, $\mathcal{E}(B)$ is dense and $\mathcal{G}_\delta$ set in the space $\bigcap_{1 \leq p < q} H^p(\mathbb{B})$.

In order to complete the proof, we consider a countable dense subset $\{w_1, w_2, w_3, \ldots\}$ of $\partial \mathbb{B}$, a decreasing sequence $\varepsilon_s$, $s = 1, 2, 3, \ldots$, of positive numbers with $\varepsilon_s \to 0$, and the balls $B(w_j, \varepsilon_s)$, centered at $w_j$ and with radious $\varepsilon_s$. By the first part of the proof, each of the sets $\mathcal{E}(B(w_j, \varepsilon_s))$ is dense and $\mathcal{G}_\delta$ set in $\bigcap_{1 \leq p < q} H^p(\mathbb{B})$. It follows from Baire's theorem that the set

$$\mathcal{Y} = \bigcap_{j=1}^{\infty} \bigcap_{s=1}^{\infty} \mathcal{E}(B(w_j, \varepsilon_s)) \text{ is dense and } \mathcal{G}_\delta \text{ in the space } \bigcap_{1 \leq p < q} H^p(\mathbb{B}).$$

We claim that the set $\mathcal{Y}$ is exactly the set of the functions $f \in \bigcap_{1 \leq p < q} H^p(\mathbb{B})$ which are totally unbounded in $\mathbb{B}$. Indeed, if $f \in \mathcal{Y}$ and $U$ is an open set with $U \cap \partial \mathbb{B} \neq \emptyset$, we may choose a point $w_{j_0} \in U \cap \partial \mathbb{B}$ and an $\varepsilon_{s_0}$ so that $B(w_{j_0}, \varepsilon_{s_0}) \subset U$. Since $\sup\{|f(z)| : z \in B(w_{j_0}, \varepsilon_{s_0}) \cap \mathbb{B}\} = \infty$, it follows that $\sup\{|f(z)| : z \in U \cap \mathbb{B}\} = \infty$. Conversely, if $f \in \bigcap_{1 \leq p < q} H^p(\mathbb{B})$ and is totally unbounded then it is obvious that $f \in \mathcal{Y}$. This completes the proof.

Next we define Hardy type spaces associated to small open subsets of the sphere $\mathbb{S} = \partial \mathbb{B}$. These are local versions of the usual Hardy spaces and the main result is that generically the functions in $\bigcap_{1 \leq p < q} H^p(\mathbb{B})$ do not belong to Hardy spaces of higher order, not even locally.

***Local Hardy spaces in the unit ball of*** $\mathbb{C}^n$***.*** Let $G \subset \mathbb{S}$ be a non-empty open set (open in $\mathbb{S}$) and $1 \leq p < \infty$. A holomorphic function $f : \mathbb{B} \to \mathbb{C}$ is said to belong to the space $H^p(\mathbb{B}, G)$ if

$$\sup_{r < 1} \int_{z \in G} |f(rz)|^p \, d\sigma(z) < \infty.$$

Now we can state the following theorem.

**Theorem 2.4.** *Let* $q \in \mathbb{R}$, $q > 1$. *Then the set*

$$\mathcal{A}_q = \left\{ g \in \bigcap_{1 \leq p < q} H^p(\mathbb{B}) : g \notin H^q(\mathbb{B}, \mathbb{S} \cap B(\zeta, \varepsilon)) \text{ for any } \zeta \in \mathbb{S} \text{ and any } \varepsilon > 0 \right\}$$

*is dense and* $\mathcal{G}_\delta$ *in the space* $\bigcap_{1 \leq p < q} H^p(\mathbb{B})$.

For the proof we will need the following lemma.

**Lemma 2.5.** *If* $\zeta \in G$ *then for the functions* $f_\zeta$ *and* $\varphi_{q,\zeta}$, *defined in Lemma 2.2, we have:*

(i) $f_\zeta \notin H^n(\mathbb{B}, G)$,

*(ii)* $\varphi_{q,\zeta} \notin H^q(\mathbb{B},G)$ *for* $1 < q < \infty$.

**Proof.** Writing
$$\int_{z\in\mathbb{S}} |f_\zeta(rz)|^n d\sigma(z) = \int_{z\in G} |f_\zeta(rz)|^n d\sigma(z) + \int_{z\in\mathbb{S}-G} |f_\zeta(rz)|^n d\sigma(z) \text{ for } r<1,$$
and taking into consideration the fact that
$$\sup_{r<1} \int_{z\in\mathbb{S}} |f_\zeta(rz)|^n d\sigma(z) = \infty,$$
we see that it suffices to show that
$$\sup_{r<1} \int_{z\in\mathbb{S}-G} |f_\zeta(rz)|^n d\sigma(z) < \infty.$$
For this, let us notice that
$$|1-\langle rz,\zeta\rangle| \geq 1 - \text{Re}(\langle rz,\zeta\rangle) = 1 - r(z\cdot\zeta).$$
Thus if $z\cdot\zeta < 0$ then $|1-\langle rz,\zeta\rangle| \geq 1$, and therefore
$$\int_{z\in(\mathbb{S}-G)\cap\{z\cdot\zeta<0\}} |f_\zeta(rz)|^n d\sigma(z) \leq \sigma(\mathbb{S}).$$
On the other hand if $z\cdot\zeta \geq 0$ then $|1-\langle rz,\zeta\rangle| \geq 1 - z\cdot\zeta$. But
$$1 - z\cdot\zeta > 0 \text{ for } z\in(\mathbb{S}-G)\cap\{z\cdot\zeta \geq 0\},$$
since $\zeta \in G$ and $z\in(\mathbb{S}-G)\cap\{z\cdot\zeta \geq 0\}$ imply that $z$ cannot be equal to $\lambda\zeta$ for any $\lambda > 0$, and therefore $z\cdot\zeta < |z\|\zeta| = 1$. By the compactness of the set $(\mathbb{S}-G)\cap\{z\cdot\zeta \geq 0\}$,
$$\alpha := \inf\{1 - z\cdot\zeta : z\in(\mathbb{S}-G)\cap\{z\cdot\zeta \geq 0\}\} > 0,$$
whence
$$\int_{z\in\mathbb{S}-G} |f_\zeta(rz)|^n d\sigma(z) \leq \frac{\sigma(\mathbb{S})}{\alpha^n} \text{ for every } r<1.$$
This proves (i). Now (ii) follows from (i), if we notice that $|\varphi_{q,\zeta}| = |f_\zeta|^{n/q}$.

***Proof of Theorem 2.4.*** Let us fix a point $w \in \mathbb{S}$ and $\delta > 0$. With $X = \{r : 0 < r < 1\}$ and $\mathcal{V} = \bigcap_{1\leq p<q} H^p(\mathbb{B})$, we consider the sublinear operator $T : \mathcal{V} \to \mathbb{C}^X$ defined as follows:
$$T(f)(r) = \left(\int_{z\in\mathbb{S}\cap B(w,\delta)} |f(rz)|^q d\sigma(z)\right)^{1/q} \text{ for } f\in\mathcal{V} \text{ and } r\in X.$$
Then, for each fixed $r \in X$, the functional $T_r : \mathcal{V} \to \mathbb{C}$, $T_r(f) = T(f)(r)$, $f\in\mathcal{V}$, is continuous. Also, by Lemma 2.5(ii), the set $\mathcal{E}(w,\delta) := \{f\in\mathcal{V} : \sup\{T(f)(r) : r\in X\} = \infty\} \neq \emptyset$. Therefore, by Lemma 2.3, the set $\mathcal{E}(w,\delta)$ is dense and $\mathcal{G}_\delta$ in the space $\bigcap_{1\leq p<q} H^p(\mathbb{B})$.

In order to complete the proof, we consider a countable dense subset $\{w_1, w_2, w_3,...\}$ of $\partial\mathbb{B}$ and a decreasing sequence $\delta_s$, $s=1,2,3,...$, of positive numbers with $\delta_s \to 0$. By the first part of the proof and Baire's theorem, the set
$$\mathcal{Y} = \bigcap_{j=1}^{\infty}\bigcap_{s=1}^{\infty} \mathcal{E}(w_j,\delta_s) \text{ is dense and } \mathcal{G}_\delta \text{ in the space } \bigcap_{1\leq p<q} H^p(\mathbb{B}).$$
We claim that $\mathcal{Y} = \mathcal{A}_q$. Indeed if $f\in\mathcal{Y}$, $\zeta\in\partial\mathbb{B}$ and $\varepsilon > 0$, we may choose $w_{j_0} \in B(\zeta,\varepsilon)$ and $\delta_{s_0}$ so that $B(w_{j_0},\delta_{s_0}) \subset B(\zeta,\varepsilon)$, and since

$$\sup_{0<r<1} \int_{z\in\mathbb{S}\cap B(w_{j_0},\delta_{s_0})} |f(rz)|^q d\sigma(z) = \infty,$$

it follows that $\sup\limits_{0<r<1} \int_{z\in\mathbb{S}\cap B(\zeta,\varepsilon)} |f(rz)|^q d\sigma(z) = \infty$, i.e., $f \notin H^q(\mathbb{B}, \mathbb{S}\cap B(\zeta,\varepsilon))$. Thus $\mathcal{Y} \subset \mathcal{A}_q$, and since it is obvious that $\mathcal{A}_q \subset \mathcal{Y}$, the proof is complete.

## 3. Hardy type spaces

First let us recall the definition of Hardy spaces in the case of bounded open sets with smooth boundary. Let $\Omega \subset \mathbb{C}^n$ be a bounded open set with $C^2$ boundary and let $\rho$ be a defining function for this set, i.e., $\rho: \mathbb{C}^n \to \mathbb{R}$ is a $C^2$ function so that $\Omega = \{\rho < 0\}$, $\partial\Omega = \{\rho = 0\}$, $\mathbb{C}^n - \overline{\Omega} = \{\rho > 0\}$ and $\nabla\rho \neq 0$ at the points of $\partial\Omega$. For $p \geq 1$, the Hardy space $H^p(\Omega)$ is defined as follows:

$$H^p(\Omega) = \left\{ f: \Omega \to \mathbb{C},\ f\ holomorphic\ in\ \Omega\ so\ that\ \|f\|_{p,\rho} := \sup_{\varepsilon > 0} \left( \int_{\{\rho=-\varepsilon\}} |f(z)|^p d\sigma_\varepsilon^\rho(z) \right)^{1/p} < \infty \right\},$$

where $d\sigma_\varepsilon^\rho$ is the Euclidean surface area measure of the hypersurface $\{z \in \mathbb{C}^n : \rho(z) = -\varepsilon\}$ (with $\varepsilon > 0$ and sufficiently small). Then $H^p(\Omega)$ is independent of the defining function $\rho$. In fact if $\lambda$ is another defining function for $\Omega$, the norms $\|f\|_{p,\rho}$ and $\|f\|_{p,\lambda}$ are equivalent. This follows from the proof of [11, Lemma 3]. Let us also observe that for compact subsets $K$ of $\Omega$,

(*) $\qquad \sup\limits_{z\in K} |f(z)| \leq A(K,\rho) \|f\|_{p,\rho}, \quad f \in H^p(\Omega),$

for some constant $A(K,\rho)$. To prove this inequality we may use the representation

$$f(z) = \int_{\partial\Omega_\varepsilon} f(\zeta) P_\varepsilon(\zeta, z) d\sigma_\varepsilon^\rho(\zeta), \quad z \in K,$$

where $P_\varepsilon(\zeta, z)$ is the Poisson kernel of $\Omega_\varepsilon := \{\rho < -\varepsilon\}$ and $\varepsilon > 0$ and sufficiently small. (Once chosen, $\varepsilon$ is fixed.) Since $P_\varepsilon(\zeta, z) \preceq [dist(z, \partial\Omega_\varepsilon)]^{-2n}$, Hölder's inequality gives

$$|f(z)| \leq \left( \int_{\zeta \in \partial\Omega_\varepsilon} |f(\zeta)|^p d\sigma_\varepsilon(\zeta) \right)^{1/p} \left( \int_{\zeta \in \partial\Omega_\varepsilon} |P_\varepsilon(\zeta,z)|^{\tilde{p}} d\sigma_\varepsilon(\zeta) \right)^{1/\tilde{p}} \preceq \frac{\|f\|_{p,\rho}}{[dist(z,\partial\Omega_\varepsilon)]^{2n}} \preceq \frac{\|f\|_{p,\rho}}{[dist(z,\partial\Omega)]^{2n}},$$

(with the point $z$ restricted to the compact set $K$) and the inequality (*) follows, if $p > 1$. (For details concerning the Poisson kernel, see [5, 11].) The case $p = 1$ is simpler.

Furthermore $H^p(\Omega)$ becomes a Banach space with the norm $\|f\|_{p,\rho}$. This follows from (*) as in the case of the unit ball. (See [13, Corollary 4.19].) From the same inequality also follows the fact that convergence in $H^p(\Omega)$ implies uniform convergence on compact subsets of $\Omega$.

As in the case of the unit ball, we define a metric in the space $\bigcap\limits_{1\leq p<q} H^p(\Omega)$, for a fixed $q > 1$, as follows. We consider a sequence $1 < p_1 < p_2 < \cdots < p_j < \cdots < q$ and $p_j \to q$, and we define the metric

$$d(f,g) = \sum_{j=1}^\infty \frac{1}{2^j} \frac{\|f-g\|_{p_j,\rho}}{1+\|f-g\|_{p_j,\rho}}.$$

Then the topology of this metric, induced on the space $\bigcap_{1\leq p<q} H^p(\Omega)$, does not depend on the choice of the sequence $p_j$ or the defining function $\rho$. Indeed, a sequence $f_k$ converges to $f$, in $\bigcap_{1\leq p<q} H^p(\mathbb{B})$, if and only if $\|f_k - f\|_{p,\rho} \to 0$, for every $p < q$.

*Local Hardy spaces.* Now we define local Hardy spaces. With $\Omega$ and $\rho$ being as above, we consider an open set $U \subset \mathbb{C}^n$ with $U \cap \partial\Omega \neq \emptyset$ and we define the space $H^p_\rho(\Omega, U)$ to be the set of holomorphic functions $f : \Omega \to \mathbb{C}$ so that $\sup_{\varepsilon > 0} \int_{\{\rho = -\varepsilon\} \cap U} |f(z)|^p d\sigma^\rho_\varepsilon(z) < \infty$. The space $H^p_\rho(\Omega, U)$ may depend on $\rho$. However we have the following lemma.

**Lemma 3.1.** *Let $\rho$ and $\lambda$ be two defining functions for $\Omega$. If $U$ and $V$ are two open subsets of $\mathbb{C}^n$ with $U \cap \partial\Omega \neq \emptyset$ and $V \cap \partial\Omega \neq \emptyset$, and if $V \subset\subset U$ then*
$$H^p_\rho(\Omega, U) \subset H^p_\lambda(\Omega, V).$$

**Proof.** The following proof is essentially the proof of [11, Lemma 3] with some minor modifications. There exist positive constants $\kappa$, $\kappa_1$ and $\kappa_2$ (independent of $\varepsilon$) so that if $z \in \{\lambda = -\varepsilon\}$ (i.e. $\lambda(z) = -\varepsilon$) then
$$B(z, \kappa\varepsilon) \subset \Lambda_\varepsilon := \{w \in \mathbb{C}^n : -\kappa_1\varepsilon < \rho(w) < -\kappa_2\varepsilon\}.$$
(The positive parameter $\varepsilon$ is assumed to be sufficiently small so that the various assertions in this proof hold true.) By the submean value property, if $f \in H^p_\rho(\Omega, U)$,
$$|f(z)|^p \leq \frac{\kappa_3}{\varepsilon^{2n}} \int_{w \in \mathbb{C}^n} \chi_\varepsilon(z, w) |f(w)|^p dw \text{ for } z \in \{\lambda = -\varepsilon\},$$
where $\chi_\varepsilon(z, w) = 1$ for $w \in B(z, \kappa\varepsilon)$ and $\chi_\varepsilon(z, w) = 0$ for $w \in \mathbb{C}^n - B(z, \kappa\varepsilon)$. ($\kappa_3$ is an appropriate constant independent of $\varepsilon$ and so are the constants $\kappa_j$, $j = 4, 5, 6$, which will be involved below.) Therefore
$$\int_{\{\lambda=-\varepsilon\}\cap V} |f(z)|^p d\sigma^\lambda_\varepsilon(z) \leq \frac{\kappa_3}{\varepsilon^{2n}} \int_{w \in \mathbb{C}^n} \left( \int_{\{\lambda=-\varepsilon\}\cap V} \chi_\varepsilon(z, w) d\sigma^\lambda_\varepsilon(z) \right) |f(w)|^p dw,$$
where we used Fubini's theorem and the measurability of the function $\chi_\varepsilon(z, w)$ for $(z, w) \in \{\lambda = -\varepsilon\} \times \mathbb{C}^n$ with respect to the product measure $d\sigma^\lambda_\varepsilon(z) \times dw$. Since $V \subset\subset U$, making $\varepsilon$ smaller – if necessary – we may assume that
$$B(z, \kappa\varepsilon) \subset \Lambda_\varepsilon \cap U \text{ for } z \in \{\lambda = -\varepsilon\} \cap V.$$
Then
$$\int_{\{\lambda=-\varepsilon\}\cap V} \chi_\varepsilon(z, w) d\sigma^\lambda_\varepsilon(z) = 0 \text{ if } w \in \mathbb{C}^n - (\Lambda_\varepsilon \cap U) \text{ and } \int_{\{\lambda=-\varepsilon\}\cap V} \chi_\varepsilon(z, w) d\sigma^\lambda_\varepsilon(z) \leq \kappa_4 \varepsilon^{2n-1} \text{ for } w \in \Lambda_\varepsilon \cap U.$$
It follows that
$$\int_{\{\lambda=-\varepsilon\}\cap V} |f(z)|^p d\sigma^\lambda_\varepsilon(z) \leq \frac{\kappa_5}{\varepsilon} \int_{w \in \Lambda_\varepsilon \cap U} |f(w)|^p dw \leq \frac{\kappa_6}{\varepsilon} \int_{\kappa_2\varepsilon}^{\kappa_1\varepsilon} \left( \int_{\{\rho=-\eta\}\cap U} |f(w)|^p d\sigma^\rho_\eta(w) \right) d\eta.$$
(The existence of the constant $\kappa_6$ follows from the coarea formula (See [2].)) Thus
$$\sup_{\varepsilon > 0} \int_{\{\lambda=-\varepsilon\}\cap V} |f(z)|^p d\sigma^\lambda_\varepsilon(z) \leq \kappa_6(\kappa_1 - \kappa_2) \sup_{\eta > 0} \int_{\{\rho=-\eta\}\cap U} |f(z)|^p d\sigma^\rho_\eta(z),$$

and this implies that $f \in H^p_\lambda(\Omega, V)$.

## 4. The case of strictly pseudoconvex domains

In this section we will show that some functions which are defined in terms of Henkin's support function belong to certain Hardy spaces. First we describe Henkin's support function $\Phi(z,\zeta)$ which is constructed in [4]. Following Henkin and Leiterer [4], let us consider an open set $\Theta \subset\subset \mathbb{C}^n$ and a $C^2$ strictly plurisubharmonic function $\rho$ in a neighbourhood of $\overline{\Theta}$. If we set

$$\beta = \frac{1}{3}\min\left\{\sum_{1\le j,k\le n}\frac{\partial^2 \rho(\zeta)}{\partial \zeta_j \partial \overline{\zeta}_k}\xi_j\overline{\xi}_k : \zeta \in \overline{\Theta},\ \xi \in \mathbb{C}^n \text{ with } |\xi|=1\right\}$$

then $\beta > 0$ and there exist $C^1$ functions $a_{jk}$ in a neighbourhood of $\overline{\Theta}$ such that

$$\max\left\{\left|a_{jk}(\zeta) - \frac{\partial^2 \rho(\zeta)}{\partial \zeta_j \partial \overline{\zeta}_k}\right| : \zeta \in \overline{\Theta}\right\} < \frac{\beta}{n^2}.$$

Let $\eta > 0$ be sufficiently small so that

$$\max\left\{\left|\frac{\partial^2 \rho(\zeta)}{\partial x_j \partial x_k} - \frac{\partial^2 \rho(z)}{\partial x_j \partial x_k}\right| : \zeta, z \in \overline{\Theta} \text{ with } |\zeta - z| \le \eta\right\} < \frac{\beta}{2n^2} \text{ for } j,k = 1,2,\ldots,2n,$$

where $x_j = x_j(\xi)$ are the real coordinates of $\xi \in \mathbb{C}^n$ such that $\xi_j = x_j(\xi) + ix_{j+n}(\xi)$. For $z, \zeta \in \overline{\Theta}$ we consider the modified Levi polynomial

$$Q(z,\zeta) = -\left[2\sum_{j=1}^n \frac{\partial \rho(\zeta)}{\partial \zeta_j}(z_j - \zeta_j) + \sum_{1\le j,k\le n} a_{jk}(z_j - \zeta_j)(z_k - \zeta_k)\right].$$

Then we have the estimate

$$\operatorname{Re} Q(z,\zeta) \ge \rho(\zeta) - \rho(z) + \beta|\zeta - z|^2 \text{ for } z, \zeta \in \overline{\Theta} \text{ with } |\zeta - z| \le \eta.$$

The following theorem is proved in [4].

**Theorem 4.1.** *Let $\Omega \subset\subset \mathbb{C}^n$ be a strictly pseudoconvex open set with $C^2$ boundary, let $\Theta$ be an open neighbourhood of $\partial\Omega$, and let $\rho$ be a $C^2$ strictly plurisubharmonic function in a neighbourhood of $\overline{\Theta}$ such that $\nabla\rho \ne 0$ at the points of $\partial\Omega$ and*

$$\Omega \cap \Theta = \{z \in \Theta : \rho(z) < 0\}.$$

*Let us choose $\eta$, $\beta$, and $Q(z,\zeta)$, as above, and let us make the positive number $\eta$ smaller so that*

$$\{z \in \mathbb{C}^n : |\zeta - z| \le 2\eta\} \subseteq \Theta \text{ for every } \zeta \in \partial\Omega.$$

*Then there exists a function $\Phi(z,\zeta)$ defined for $\zeta$ in some open neighbourhood $U_{\partial\Omega} \subseteq \Theta$ of $\partial\Omega$ and $z \in U_{\overline{\Omega}} = \Omega \cup U_{\partial\Omega}$, which is $C^1$ in $(z,\zeta) \in U_{\overline{\Omega}} \times U_{\partial\Omega}$, holomorphic in $z \in U_{\overline{\Omega}}$, and such that $\Phi(z,\zeta) \ne 0$ for $(z,\zeta) \in U_{\overline{\Omega}} \times U_{\partial\Omega}$ with $|\zeta - z| \ge \eta$, and*

$$\Phi(z,\zeta) = Q(z,\zeta)C(z,\zeta) \text{ for } (z,\zeta) \in U_{\overline{\Omega}} \times U_{\partial\Omega} \text{ with } |\zeta - z| \le \eta,$$

*for some $C^1$-function $C(z,\zeta)$ defined for $(z,\zeta) \in U_{\overline{\Omega}} \times U_{\partial\Omega}$ and $\ne 0$ when $|\zeta - z| \le \eta$.*

In this setting we will prove the following lemma. We use a set of coordinates – *the Levi coordinates* – which are appropriate when we are dealing with integrals involving the function $\Phi(z,\zeta)$. (See [4, 8, 12] .) As a matter of fact we will use a slight modification of the Levi coordinates.

**Lemma 4.2.** *For each fixed point* $\zeta \in \partial\Omega$,

(1) $$\sup_{\varepsilon > 0} \int_{\{\rho(z)=-\varepsilon\}} \frac{d\sigma_\varepsilon(z)}{|\Phi(z,\zeta)|^p} < \infty \text{ when } 1 < p < n, \text{ and}$$

(2) $$\sup_{\varepsilon > 0} \int_{\{\rho(z)=-\varepsilon\}} \frac{d\sigma_\varepsilon(z)}{|\Phi(z,\zeta)|^{2n-1}} = \infty.$$

*Therefore* $\dfrac{1}{\Phi(\cdot,\zeta)} \in \bigcap_{1 \le p < n} H^p(\Omega)$ *and* $\dfrac{1}{\Phi(\cdot,\zeta)} \notin H^{2n-1}(\Omega)$.

**Proof.** We consider a coordinate system $t = (t_1, t_2, t_3, \ldots, t_{2n}) = (t_1(z), t_2(z), t_3(z), \ldots, t_{2n}(z))$, of real $C^1$–functions, for points $z \in \mathbb{C}^n = \mathbb{R}^{2n}$, which are sufficiently close to the point $\zeta$, as follows: We set

$$t_1(z) = -\rho(z) \text{ and } t_2(z) = \text{Im} Q(z,\zeta).$$

Then $d_z Q(z,\zeta)\big|_{z=\zeta} = -2 \sum_{j=1}^{n} \dfrac{\partial \rho(\zeta)}{\partial \zeta_j} dz_j \bigg|_{z=\zeta} = -2\partial\rho(\zeta)$ and, therefore,

$$d_z t_2(z)\big|_{z=\zeta} = d_z[\text{Im} Q(z,\zeta)]\big|_{z=\zeta} = i[\partial\rho(\zeta) - \bar\partial\rho(\zeta)].$$

On the other hand,

$$d_z t_1(z)\big|_{z=\zeta} = d_z[-\rho(z)]\big|_{z=\zeta} = -[\partial\rho(\zeta) + \bar\partial\rho(\zeta)].$$

It follows that

$$\left(d_z t_1(z)\big|_{z=\zeta}\right) \wedge \left(d_z t_2(z)\big|_{z=\zeta}\right) = 2i\partial\rho(\zeta) \wedge \bar\partial\rho(\zeta) \ne 0.$$

Now the existence of $C^1$–functions $t_3(z), \ldots, t_{2n}(z)$ such that the mapping

$$z \to (t_1(z), t_2(z), t_3(z), \ldots, t_{2n}(z))$$

is a $C^1$–diffeomorphism, from an open neighbourhood of the point $\zeta$ to an open neighbourhood of $0 \in \mathbb{C}^n = \mathbb{R}^{2n}$ (with $t(\zeta) = 0$), follows from the inverse function theorem. Also let us point out that, for $z$ sufficiently close to $\zeta$, $z \in \Omega$ if and only if $t_1 = -\rho(z) > 0$. For points $z \in \Omega$ which are sufficiently close to $\zeta$,

$$|\Phi(z,\zeta)| \approx |Q(z,\zeta)| \approx |\text{Re}\, Q(z,\zeta)| + |\text{Im}\, Q(z,\zeta)| \ge -\rho(z) + \beta|\zeta - z|^2 + |\text{Im}\, Q(z,\zeta)|$$

and

$$|\zeta - z|^2 \approx t_1^2 + t_2^2 + t_3^2 + \cdots + t_{2n}^2.$$

(When we write $A \approx B$, we mean that $\kappa B \le A \le \mu B$, for some positive constants $\kappa$ and $\mu$ which are independent of $z$.) Therefore (for $z \in \Omega$ and sufficiently close to $\zeta$)

$$|\Phi(z,\zeta)| \succeq t_1 + t_1^2 + t_2^2 + t_3^2 + \cdots + t_{2n}^2 + |t_2|.$$

Therefore (1) follows from

$$\sup_{\varepsilon > 0} \int_{t_1 = \varepsilon} \frac{dt_2 \cdots dt_{2n}}{(t_1 + |t_2| + t_1^2 + t_2^2 + t_3^2 + \cdots + t_{2n}^2)^p} < \infty$$

or equivalently from

(3) $$\sup_{\varepsilon > 0} \int_{t_1 = \varepsilon} \frac{dt_2 \cdots dt_{2n}}{(t_1 + |t_2| + t_3^2 + \cdots + t_{2n}^2)^p} < \infty.$$

(In the above integrals $t$ is restricted in a «small» neighbourhood of $0 \in \mathbb{C}^n = \mathbb{R}^{2n}$, i.e., $|t|$ is «small».) But

$$\int\limits_{t_1=\varepsilon} \frac{dt_2 \cdots dt_{2n}}{(t_1+|t_2|+t_3^2+\cdots+t_{2n}^2)^p} \approx \int\limits_{t_2^2+t_3^2+\cdots+t_{2n}^2<1} \frac{dt_2 \cdots dt_{2n}}{(\varepsilon+|t_2|+t_3^2+\cdots+t_{2n}^2)^p}.$$

Also, by Fubini's theorem,

$$\int\limits_{t_2^2+t_3^2+\cdots+t_{2n}^2<1} \frac{dt_2 \cdots dt_{2n}}{(\varepsilon+|t_2|+t_3^2+\cdots+t_{2n}^2)^p} \approx \int\limits_{t_3^2+\cdots+t_{2n}^2<1} \left( \int\limits_{t_2=0}^{1} \frac{dt_2}{(\varepsilon+t_2+t_3^2+\cdots+t_{2n}^2)^p} \right) dt_3 \cdots dt_{2n}$$

$$\approx \int\limits_{t_3^2+\cdots+t_{2n}^2<1} \frac{dt_3 \cdots dt_{2n}}{(\varepsilon+t_3^2+\cdots+t_{2n}^2)^{p-1}}.$$

Integrating in polar coordinates we see that the last integral is equal to

$$\int\limits_{r=0}^{1} \frac{r^{2n-3} dr}{(\varepsilon+r^2)^{p-1}} \leq \int\limits_{r=0}^{1} r^{2n-2p-1} dr < \infty.$$

This proves (3) and completes the proof of (1).

In order to prove (2), let us observe that for points $z \in \Omega$ which are sufficiently close to $\zeta$,

$$|\Phi(z,\zeta)| \approx |Q(z,\zeta)| \preceq |\zeta - z| \approx (t_1^2+t_2^2+t_3^2+\cdots+t_{2n}^2)^{1/2},$$

whence

$$\int\limits_{\{\rho(z)=-\varepsilon\}} \frac{d\sigma_\varepsilon(z)}{|\Phi(z,\zeta)|^{2n-1}} \approx \int\limits_{t_1=\varepsilon} \frac{dt_2 \cdots dt_{2n}}{(t_1^2+t_2^2+\cdots+t_{2n}^2)^{(2n-1)/2}} \succeq \int\limits_{t_2^2+t_3^2+\cdots+t_{2n}^2<1} \frac{dt_2 \cdots dt_{2n}}{(\varepsilon^2+t_2^2+\cdots+t_{2n}^2)^{(2n-1)/2}}.$$

By introducing polar coordinates in the last integral, we see that this integral is equal to

$$\text{constant} \int\limits_{r=0}^{1} \frac{r^{2n-2} dr}{(\varepsilon^2+r^2)^{(2n-1)/2}} \approx \int\limits_{r=0}^{1} \frac{r^{2n-2} dr}{(\varepsilon+r)^{2n-1}}.$$

But as $\varepsilon$ decreases, the function $r^{2n-2}/(\varepsilon+r)^{2n-1}$ increases, and the monotone convergence theorem gives that

$$\lim_{\varepsilon \to 0^+} \int\limits_{r=0}^{1} \frac{r^{2n-2} dr}{(\varepsilon+r)^{2n-1}} = \int\limits_{r=0}^{1} \frac{dr}{r} = \infty,$$

and proves (2).

**Lemma 4.3.** *Let $\Omega \subset\subset \mathbb{C}^n$ be a strictly pseudoconvex open set with $C^2$ boundary and let $\rho$ be a $C^2$ strictly plurisubharmonic defining function of $\Omega$ defined in a neighbourhood of $\overline{\Omega}$. If $1 < q < \infty$ and $U \subset \mathbb{C}^n$ with $U \cap \partial \Omega \neq \varnothing$, then there exists a function $h_{q,U}$ so that*

$$h_{q,U} \in \bigcap_{1 \leq p < q} H^p(\Omega) \text{ and } h_{q,U} \notin H_\rho^{(2n-1)q/n}(\Omega, U).$$

**Proof.** Let us fix a point $\zeta \in U \cap \partial\Omega$. Then, as it follows from Taylor's theorem and the strict plurisubharmonicity of $\rho$ (see [8]), the Levi polynomial of $\rho$

$$F(z,\zeta) = -\left[ 2\sum_{j=1}^{n} \frac{\partial \rho(\zeta)}{\partial \zeta_j}(z_j-\zeta_j) + \sum_{1 \leq j,k \leq n} \frac{\partial^2 \rho(\zeta)}{\partial \zeta_j \partial \zeta_k}(z_j-\zeta_j)(z_k-\zeta_k) \right]$$

satisfies the inequality

$$\text{Re}\, F(z,\zeta) \geq \rho(\zeta) - \rho(z) + \gamma |\zeta-z|^2 \text{ for } z \in \mathbb{C}^n \text{ with } |\zeta-z| < \delta,$$

for some «small» positive constants $\delta$ and $\gamma$. In particular,

$$\text{Re}\, F(z,\zeta) > 0 \text{ for } z \in B(\zeta, \delta) \cap \overline{\Omega} - \{\zeta\}.$$

It follows that the function $\log[1/F(z,\zeta)]$ is defined and holomorphic for $z \in B(\zeta,\delta) \cap \Omega$, and that $\lim_{z \in \Omega, z \to \zeta} \log[1/F(z,\zeta)] = \infty$. (Here log is the principal branch of the logarithm with $|\arg| \leq \pi$.) Also we can prove, as in the proof of the Lemma 4.2, that if $q < n$,

(1) $$\sup_{\varepsilon > 0} \int_{\{\rho(z) = -\varepsilon\} \cap B(\zeta,\delta)} \frac{d\sigma_\varepsilon(z)}{|F(z,\zeta)|^q} < \infty \text{ and } \sup_{\varepsilon > 0} \int_{\{\rho(z) = -\varepsilon\} \cap B(\zeta,\delta)} \frac{d\sigma_\varepsilon(z)}{|F(z,\zeta)|^{2n-1}} = \infty.$$

Then, using the first of part of (1), we obtain, as in Lemma,

(2) $$\sup_{\varepsilon > 0} \int_{\{\rho(z) = -\varepsilon\} \cap B(\zeta, 2\delta/3)} \left|\log\left[\frac{1}{F(z,\zeta)}\right]\right|^p d\sigma_\varepsilon(z) < \infty, \text{ for every } p < \infty.$$

Next we consider a $C^\infty$-function $\chi : \mathbb{C}^n \to \mathbb{R}$, $0 \leq \chi(z) \leq 1$, with compact support contained in $B(\zeta, 2\delta/3)$, and such that $\chi(z) = 1$ when $z \in B(\zeta, \delta/3)$. Now the function

$$\chi(z) \log\left[\frac{1}{F(z,\zeta)}\right]$$

is extended to a $C^\infty$-function in $\Omega$, by defining it to be 0 in $\Omega - B(\zeta, 2\delta/3)$. Then the $(0,1)$-form

$$u(z) := \bar{\partial}\left\{\chi(z) \log\left[\frac{1}{F(z,\zeta)}\right]\right\}$$

is defined and is $C^\infty$ in a open neighbourhood $\bar{\Omega}$, it is zero for $z \in B(\zeta, \delta/3) \cap \Omega$, and, in particular, it has bounded coefficients in $\Omega$. In fact $u(z)$ extends to a $C^\infty$ $(0,1)$-form for $z$ in an open neighbourhood of $\bar{\Omega}$, since the function $\log\left[\frac{1}{F(z,\zeta)}\right]$ is holomorphic in an open neighbourhood of the compact set $\overline{[B(\zeta, 2\delta/3) - B(\zeta, \delta/3)] \cap \Omega}$. It follows that there exists a bounded $C^\infty$-function $\psi : \Omega \to \mathbb{C}$ which solves the equation $\bar{\partial}\psi = u$ in $\Omega$. (See [8].) Then we may define the functions

$$f_\zeta(z) := \chi(z) \log\left[\frac{1}{F(z,\zeta)}\right] - \psi(z)$$

and

$$h_{q,\zeta}(z) := \exp\left[\frac{n}{q} f_\zeta(z)\right] = \exp\left\{\chi(z) \log\left[\frac{1}{F(z,\zeta)}\right]^{n/q} - \frac{n}{q}\psi(z)\right\}.$$

Then the functions $f_\zeta(z)$ and $h_{q,\zeta}(z)$ are holomorphic for $z \in \Omega$. We claim that

(1) $$\sup_{\varepsilon > 0} \int_{\{\rho(z) = -\varepsilon\} \cap B(\zeta,\delta)} |h_{q,\zeta}(z)|^p d\sigma_\varepsilon^\rho(z) < \infty \text{ for } p < q.$$

Notice that the behaviour of the above integral is not affected by the functions $\chi$ or $\psi$, since $\chi \equiv 1$ near $\zeta$ and $\psi$ is bounded in $\Omega$. Thus (2) follows from the first part of the (1).
Also

$$\sup_{\varepsilon > 0} \int_{\{\rho(z) = -\varepsilon\} \cap B(\zeta,\delta)} |h_{q,\zeta}(z)|^{(2n-1)q/n} d\sigma_\varepsilon^\rho(z) = \infty.$$

Indeed this follows from the second part of (1), since $\chi \equiv 1$ near $\zeta$ and $\exp(-\psi)$ is bounded away from zero in $\Omega$.
Thus setting $h_{q,U} := h_{q,\zeta}$ we obtain the required function.

**Theorem 4.4.** *Let* $\Omega \subset\subset \mathbb{C}^n$ *be a strictly pseudoconvex open set with* $C^2$ *boundary and* $q \in \mathbb{R} \cup \{\infty\}$, $q > 1$. *Then the following hold:*

*(i) The set of the functions in the space $\bigcap_{1\leq p<q} H^p(\Omega)$ which are totally unbounded in $\Omega$ is dense and $\mathcal{G}_\delta$ in this space.*

*(ii) The set of the functions in the space $\bigcap_{1\leq p<q} H^p(\Omega)$ which are singular at every boundary point of $\Omega$ is dense and $\mathcal{G}_\delta$ in this space.*

**Proof.** Let us consider a 'small' ball $B$ whose center lies on $\partial\Omega$, and let us set $X = B \cap \Omega$ and $\mathcal{V} = \bigcap_{1\leq p<q} H^p(\Omega)$. We define the linear operator

$$T : \mathcal{V} \to \mathbb{C}^X \text{ with } T(f)(z) = f(z) \text{ for } z \in X \text{ and } f \in \mathcal{V}.$$

For each fixed $z \in X$, the functional $T_z : \mathcal{V} \to \mathbb{C}$ defined by $T_z(f) = f(z)$, $f \in \mathcal{V}$, is continuous. It is easy to see that the set $\mathcal{E} = \{f \in \mathcal{V} : T(f) \text{ is unbounded on } X\}$ in this case is equal to

$$\mathcal{E}(B) = \left\{ f \in \bigcap_{1\leq p<q} H^p(\Omega) : \sup_{z \in B \cap \Omega} |f(z)| = \infty \right\}.$$

Now we consider the function $f_\zeta$ which was constructed in the proof of Lemma 4.3. If the point $\zeta \in B \cap \partial\Omega$ then $f_\zeta \in \mathcal{E}(B)$, and therefore $\mathcal{E}(B) \neq \varnothing$. Therefore, by Lemma 2.3, $\mathcal{E}(B)$ is dense and $\mathcal{G}_\delta$ set in the space $\mathcal{V}$. In order to complete the proof, we consider a countable dense subset $\{w_1, w_2, w_3, ...\}$ of $\partial\Omega$, a decreasing sequence $\varepsilon_s$, $s = 1,2,3,...$, of positive numbers with $\varepsilon_s \to 0$, and the balls $B(w_j, \varepsilon_s)$, centered at $w_j$ and with radious $\varepsilon_s$. By the first part of the proof, each of the sets $\mathcal{E}(B(w_j, \varepsilon_s))$ is dense and $\mathcal{G}_\delta$ set in $\bigcap_{1\leq p<q} H^p(\Omega)$. It follows from Baire's theorem that the set

$$\mathcal{Y} = \bigcap_{j=1}^{\infty} \bigcap_{s=1}^{\infty} \mathcal{E}(B(w_j, \varepsilon_s)) \text{ is dense and } \mathcal{G}_\delta \text{ in the space } \bigcap_{1\leq p<q} H^p(\Omega).$$

It is easy to see that $\mathcal{Y}$ is the set of the functions in the space $\bigcap_{1\leq p<q} H^p(\Omega)$ which are totally unbounded in $\Omega$, and this proves (i). The assertion (ii) follows from [6].

**Theorem 4.5.** *Let $\Omega \subset\subset \mathbb{C}^n$ be a strictly pseudoconvex open set with $C^2$ boundary and let $\rho$ be a $C^2$ strictly plurisubharmonic defining function of $\Omega$ defined in a neighbourhood of $\overline{\Omega}$. If $q \in \mathbb{R}$, $q > 1$, then the set*

$$\mathcal{B}_q^\rho = \left\{ g \in \bigcap_{1\leq p<q} H^p(\Omega) : g \notin H_\rho^{(2n-1)q/n}(\Omega, U) \text{ for any open set } U \text{ with } U \cap \partial\Omega \neq \varnothing \right\}$$

*is dense and $\mathcal{G}_\delta$ in the space $\bigcap_{1\leq p<q} H^p(\Omega)$.*

**Proof.** Let us fix a point $w \in \partial\Omega$ and a positive number $\delta$. With $X = \{\varepsilon : 0 < \varepsilon < \varepsilon_0\}$ (where $\varepsilon_0$ is a 'small' positive number) and $\mathcal{V} = \bigcap_{1\leq p<q} H^p(\Omega)$, we consider the sublinear operator $T : \mathcal{V} \to \mathbb{C}^X$ defined as follows:

$$T(f)(\varepsilon) = \left( \int_{\{\rho=-\varepsilon\} \cap B(w,\delta)} |f(z)|^{(2n-1)q/n} d\sigma_\varepsilon^\rho(z) \right)^{\frac{n}{(2n-1)q}} \text{ for } f \in \mathcal{V} \text{ and } \varepsilon \in X.$$

Then, for each fixed $\varepsilon \in X$, the functional $T_\varepsilon : \mathcal{V} \to \mathbb{C}$, $T_\varepsilon(f) = T(f)(\varepsilon)$, $f \in \mathcal{V}$, is continuous. Also, by Lemma 4.3, the set $\mathcal{E}(w,\delta) := \{f \in \mathcal{V} : \sup\{T(f)(\varepsilon) : \varepsilon \in X\} = \infty\} \neq \varnothing$. Therefore, by Lemma 2.3, the set $\mathcal{E}(w,\delta)$ is dense and $\mathcal{G}_\delta$ in the space $\mathcal{V}$.

In order to complete the proof, we consider a countable dense subset $\{w_1, w_2, w_3, ...\}$ of $\partial\Omega$ and a decreasing sequence $\delta_s$, $s = 1,2,3,...$, of positive numbers with $\delta_s \to 0$. By the first part of the proof and Baire's theorem, the set

$$\mathcal{Y} = \bigcap_{j=1}^{\infty} \bigcap_{s=1}^{\infty} \mathcal{E}(w_j, \delta_s) \text{ is dense and } \mathcal{G}_\delta \text{ in the space } \bigcap_{1 \leq p < q} H^p(\Omega).$$

Now it easy to see that $\mathcal{Y} = \mathcal{B}_q^\rho$, and this completes the proof.

Combining Theorem 4.5 with Lemma 3.1, we see that the set $\mathcal{B}_q^\rho$ is independant of $\rho$. Thus we have the following theorem.

**Theorem 4.6.** *Let $\Omega \subset\subset \mathbb{C}^n$ be a strictly pseudoconvex open set with $C^2$ boundary. If $q \in \mathbb{R}$, $q > 1$, then the set*

$$\mathcal{B}_q = \{g \in \bigcap_{1 \leq p < q} H^p(\Omega) : g \notin H_\lambda^{(2n-1)q/n}(\Omega, U) \text{ for any open set } U \text{ with } U \cap \partial\Omega \neq \varnothing$$

*and any defining function $\lambda$ of $\Omega\}$*

*is dense and $\mathcal{G}_\delta$ in the space $\bigcap_{1 \leq p < q} H^p(\Omega)$.*

**Remark 4.7.** It is easy to see that one can prove theorems analogous to the Theorems 2.1, 2.4, 4.4, 4.6, also with the spaces $H^p$ in place of the intersections $\bigcap_{p < q} H^p$. Thus we have

1. For $1 \leq p < \infty$, the set of the functions in the space $H^p(\mathbb{B})$ which are totally unbounded in $\mathbb{B}$ is dense and $\mathcal{G}_\delta$ in this space.

2. For $1 \leq p < q < \infty$, the set

$$\{g \in H^p(\mathbb{B}) : g \notin H^q(\mathbb{B}, \mathbb{S} \cap B(\zeta, \varepsilon)) \text{ for any } \zeta \in \mathbb{S} \text{ and any } \varepsilon > 0\}$$

is dense and $\mathcal{G}_\delta$ in the space $H^p(\mathbb{B})$.

3. If $\Omega \subset\subset \mathbb{C}^n$ is a strictly pseudoconvex open set with $C^2$ boundary and $1 \leq p < \infty$, the set of the functions in the space $H^p(\Omega)$ which are totally unbounded in $\Omega$ is dense and $\mathcal{G}_\delta$ in this space.

4. If $\Omega \subset\subset \mathbb{C}^n$ is a strictly pseudoconvex open set with $C^2$ boundary, $1 \leq p < \infty$ and $q > (2n-1)p/n$, then the set

$$\{g \in H^p(\Omega) : g \notin H_\lambda^q(\Omega, U) \text{ for any open set } U \text{ with } U \cap \partial\Omega \neq \varnothing$$

*and any defining function $\lambda$ of $\Omega\}$*

is dense and $\mathcal{G}_\delta$ in the space $H^p(\Omega)$.

## 5. Hardy spaces of harmonic functions

Let $\Omega \subset \mathbb{R}^n$ be a bounded open set with $C^2$ boundary. If $\rho$ is a $C^2$ defining function of $\Omega$ then one can define the harmonic Hardy spaces $h^p(\Omega)$, $p \geq 1$, (see [1, 11]), the intersections $\bigcap_{p < q} h^p(\Omega)$, and the local Hardy spaces $h_\rho^p(\Omega, U)$, as before. ($U \subset \mathbb{R}^n$ is an open set with $U \cap \partial\Omega \neq \varnothing$.)

**Lemma 5.1.** Let $n \geq 3$ and $y \in \partial\Omega$. Then the function $\varphi_y(x) = \dfrac{1}{|x-y|^{n-2}}$ ($x \neq y$) belongs to $h^p(\Omega)$ if and only if $p < \dfrac{n-1}{n-2}$. In particular $\varphi_y \notin h^{(n-1)/(n-2)}(\Omega)$.

**Proof.** We may assume that $y = 0 \in \partial\Omega$, in which case $\varphi_y$ becomes the function
$$\varphi_0(x) = \frac{1}{|x|^{n-2}} = \frac{1}{(x_1^2 + x_2^2 + \cdots + x_n^2)^{(n-2)/2}}.$$

We must show that

(1) $\quad \sup\limits_{\varepsilon > 0} \int\limits_{\{\rho = -\varepsilon\}} \dfrac{d\sigma_\varepsilon^\rho(x)}{|x|^{(n-2)p}} < \infty$ if and only if $p < \dfrac{n-1}{n-2}$.

Using a local diffeomorphism – near the point 0 of $\partial\Omega$ – we may assume that the hypersurface $\partial\Omega$, near 0, is defined by the equation $x_1 = 0$, and that $x_1 > 0$ for $x \in \Omega$ (close to 0). Then (1) is equivalent to

(2) $\quad \sup\limits_{\varepsilon > 0} \int\limits_{x_2^2 + \cdots + x_n^2 < 1} \dfrac{dx_2 \ldots dx_n}{(\varepsilon^2 + x_2^2 + \cdots + x_n^2)^{(n-2)p/2}} < \infty$ if and only if $p < \dfrac{n-1}{n-2}$.

Integrating in polar coordinates we see that the above integral becomes
$$(constant) \int_{r=0}^{1} \frac{r^{n-2} dr}{(\varepsilon^2 + r^2)^{(n-2)p/2}}.$$

By monotone convergence theorem,
$$\sup_{\varepsilon > 0} \int_{r=0}^{1} \frac{r^{n-2} dr}{(\varepsilon^2 + r^2)^{(n-2)p/2}} = \lim_{\varepsilon \to 0^+} \int_{r=0}^{1} \frac{r^{n-2} dr}{(\varepsilon^2 + r^2)^{(n-2)p/2}} \leq \int_{r=0}^{1} \frac{r^{n-2} dr}{r^{(n-2)p}},$$

and (2) follows.

**Lemma 5.2.** Let $n \geq 3$ and $y \in \partial\Omega$. Then $\varphi_y \notin h_\rho^{(n-1)/(n-2)}(\Omega, U)$ for $y \in U$.

**Proof.** It follows easily from the previous lemma.

With the above lemmas, we can prove the following theorems. Their proofs are similar to the proofs of Theorems 4.4 and 4.6.

**Theorem 5.3.** *Let $q \in \mathbb{R} \cup \{\infty\}$ with $1 < q \leq \dfrac{n-1}{n-2}$. Then the set of the functions in the space $\bigcap\limits_{1 \leq p < q} h^p(\Omega)$ which are totally unbounded in $\Omega$ is dense and $\mathcal{G}_\delta$ in this space.*

**Theorem 5.4.** *Let $q \in \mathbb{R}$ with $1 < q \leq \dfrac{n-1}{n-2}$. Then the set*
$$\mathcal{A}_q = \{g \in \bigcap_{1 \leq p < q} h^p(\Omega) : g \notin h_\lambda^{(n-1)/n-2)}(\Omega, U) \text{ for any open set } U \text{ with } U \cap \partial\Omega \neq \emptyset$$

*and any defining function $\lambda$ of $\Omega$*$\}$

*is dense and $\mathcal{G}_\delta$ in the space $\bigcap\limits_{1 \leq p < q} h^p(\Omega)$.*

**Remark 5.5.** According to Theorem 5.3, the functions in the space $\bigcap_{1\leq p<q} h^p(\Omega)$ are generically totally unbounded in $\Omega$, despite the fact that all these functions have non-tangential limits almost everywhere at the points of the boundary of $\Omega$ (by Fatou's theorem). Similar remarks can be made for Theorems 2.1 and 4.4.

**Acknowledgment.** I would like to thank T. Hatziafratis and V. Nestoridis for helpful discussions.

K. Kioulafa
National and Kapodistrian University of Athens
Department of Mathematics
Panepistemiopolis
157 84 Athens
Greece

E-mail address: keiranna@math.uoa.gr